\documentclass[12pt, dvips, twoside]{amsart}
\usepackage{amssymb,latexsym,bbm,graphicx,epsfig,epic,eepic,oldgerm}
\usepackage{psfrag}

\theoremstyle{plain}
\newtheorem{theorem}{Theorem}

\newtheorem{proposition}[theorem]{Proposition}
\newtheorem{corollary}[theorem]{Corollary}
\newtheorem{remark}[theorem]{Remark}

\newcommand{\proofend}{\hspace*{\fill} $\Box$\\}

\def\s{\smallskip}
\def\m{\medskip}
\def\b{\bigskip}
\def\ni{\noindent}

\def\dd{\delta}
\def\eps{\epsilon}
\def\ff{\varphi}

\def\ss{\sigma}

\def\RR{\mathbbm{R}}

\def\pp{\partial}

\def\ra{\rightarrow}
\def\ha{\hookrightarrow}

\def\Vol{\operatorname{Vol}\;\!}

\begin{document}

\begin{titlepage}
\title{Volume preserving embeddings of open subsets of $\RR^n$ into manifolds}    
$   $ \\
$   $ \\

\author{Felix Schlenk}
\address{(F.\ Schlenk) ETH Z\"urich, CH-8092 Z\"urich, Switzerland}
\email{felix@math.ethz.ch}

\date{\today}
\end{titlepage}

\begin{abstract}
We consider a connected smooth $n$-dimensional manifold $M$ endowed with 
a volume form $\Omega$, and we show that an open subset $U$ of $\RR^n$
of Lebesgue measure $\Vol (U)$ embeds into $M$ by a smooth volume
preserving embedding whenever the volume condition $\Vol (U) \le \Vol (M,
\Omega)$ is met.
\end{abstract}

\maketitle

\markboth{{\rm Volume preserving embeddings of open subsets of $\RR^n$
               into manifolds}}{{}} 

\section{Introduction}
\ni
Consider a connected smooth $n$-dimensional manifold $M$ with or without 
boundary.
A volume form on $M$ is a smooth nowhere vanishing differential $n$-form 
$\Omega$.
It follows that $M$ is orientable. We orient $M$ such that $\int_M
\Omega$ is positive, and we write $\Vol (M, \Omega) = \int_M \Omega$.
We endow each open (not necessarily connected) subset
$U$ of $\RR^n$ with the Euclidean volume form
\[
\Omega_0 = dx_1 \wedge \dots \wedge dx_n .
\]
A smooth embedding $\ff \colon U \ha M$ is called volume preserving if
\[
\ff^* \Omega = \Omega_0 .
\]
Then $\Vol (U, \Omega_0) \le \Vol (M, \Omega)$.
In this note we prove that this obvious condition for the existence of a volume
preserving embedding is the only one.

\begin{theorem}  \label{ta:1}
Consider an open subset $U$ of $\RR^n$ and a smooth connected 
$n$-dimensional manifold $M$ endowed with a volume form 
$\Omega$. 
Then there exists a volume preserving embedding $\ff \colon U \ha M$ if 
and only if  
$\Vol (U, \Omega_0) \le \Vol (M, \Omega).$ 
\end{theorem}

\ni
If $U$ is a bounded subset whose boundary has zero measure and if 
$\Vol (U, \Omega_0) < \Vol (M, \Omega)$,
Theorem \ref{ta:1} is an easy consequence of 
Moser's deformation method.
Moreover, if $U$ is a ball and $M$ is compact, Theorem \ref{ta:1} has been
proved in \cite{K}.
The main point of this note therefore is to show that Theorem \ref{ta:1} 
holds true for an arbitrary open subset of $\RR^n$ and an arbitrary connected
manifold even in case that the volumes are equal.

\b
\ni
{\it Acknowledgement.}
I'm grateful to Tom Ilmanen and Edi Zehnder for helpful discussions.

\setcounter{secnumdepth}{4}
\setcounter{equation}{0}

\section{Proof of Theorem \ref{ta:1}}

\noindent
Assume first that $\ff \colon U \ha M$ is a smooth embedding such that $\ff^* \Omega = \Omega_0$. Then
\[
\Vol (U, \Omega_0) \,=\, \int_U \Omega_0 \,=\, \int_U \ff^* \Omega \,=\, \int_{\ff(U)} \Omega \,\le\, \int_M \Omega \,=\, \Vol (M, \Omega) .
\]

\s
\ni
Assume now that $\Vol (U, \Omega_0) \,\le\,  \Vol (M, \Omega)$.
We are going to construct a smooth embedding $\ff \colon U \ha M$ such
that $\ff^* \Omega = \Omega_0$.

We orient $\RR^n$ in the natural way.
The orientations of $\RR^n$ and $M$ orient each open subset of $\RR^n$ and $M$.
We abbreviate the Lebesgue measure $\Vol (V, \Omega_0)$ of a measurable
subset $V$ of $\RR^n$ by $\left| V \right|$,
and we write $\overline{V}$ for the closure of $V$ in $\RR^n$.
Moreover, we denote by $B_r$ the open ball in $\RR^n$ of radius $r$
centered at the origin.

\begin{proposition}  \label{l:smooth}
Assume that $V$ is a non-empty open subset of $\RR^n$. Then there exists a smooth embedding $\ss \colon V \ha \RR^n$ such that 
$\left| \RR^n \setminus \ss (V) \right| =0$.
\end{proposition}

\proof
We choose an increasing sequence
\[
V_1 \subset V_2 \subset \cdots \subset V_k \subset V_{k+1} \subset \cdots
\] 
of non-empty open subsets of $V$ such that $\overline{V_k} \subset V_{k+1}$, 
$k=1,2, \dots$, and $\bigcup_{k=1}^\infty V_k = V$.
To fix the ideas, we assume that the sets $V_k$ have smooth boundaries.

Let $\ss_1 \colon V_2 \ha \RR^n$ be a smooth embedding such that 
$\ss_1(V_1) \subset B_1$ and
\[
\left| B_1 \setminus \ss_1(V_1) \right| \,\le\, 2^{-1} .
\]
Since $\overline{V_1} \subset V_2$ and $\overline{\ss_1(V_1)} \subset \overline{B_1} \subset B_2$, we find a smooth embedding $\ss_2 \colon V_3 \ha \RR^n$
such that 
$\ss_2 |_{V_1} = \ss_1 |_{V_1}$ and $\ss_2(V_2) \subset B_2$ and
\[
\left| B_2 \setminus \ss_2(V_2) \right| \,\le\, 2^{-2} .
\]
Arguing by induction we find smooth embeddings $\ss_k \colon V_{k+1} \ha \RR^n$
such that 
$\ss_k |_{V_{k-1}} = \ss_{k-1} |_{V_{k-1}}$ and $\ss_k(V_k) \subset B_k$ and
\begin{equation}  \label{e:-k}
\left| B_k \setminus \ss_k(V_k) \right| \,\le\, 2^{-k} ,
\end{equation}
$k =1, 2, \dots$.
The map $\ss \colon V \ra \RR^n$ defined by 
$ \ss |_{V_k} = \ss_k |_{V_k}$
is a well defined smooth embedding of $V$ into $\RR^n$.
Moreover, the inclusions $\ss_k(V_k) \subset \ss(V)$ and the estimates \eqref{e:-k} imply that
\[
\left| B_k \setminus \ss(V) \right| \,\le\, \left| B_k \setminus \ss_k(V_k) \right| \,\le\, 2^{-k},
\]
and so 
\[
\left| \RR^n \setminus \ss (V) \right| \,=\, \lim_{k \ra \infty} \left| B_k \setminus \ss(V) \right| \,=\, 0 .
\]
This completes the proof of Proposition \ref{l:smooth}.
\proofend

\ni
Our next goal is to construct a smooth embedding of $\RR^n$ into the
connected $n$-dimensional manifold $M$ such that the complement of the
image has measure zero. 
If $M$ is compact, such an embedding has been obtained by Ozols \cite{O} 
and Katok \cite[Proposition 1.3]{K}.
While Ozols combines an engulfing method with tools from Riemannian
geometry, Katok successively exhausts a smooth triangulation of $M$.
Both approaches can be generalized to the case of an arbitrary connected 
manifold $M$, and we shall follow Ozols.

\s
We abbreviate $\RR_{>0} = \{\:\! r \in \RR \mid r >0 \:\!\}$ and 
$\overline{\RR}_{>0} = \RR_{>0}\cup \{ \infty \}$.
We endow $\overline{\RR}_{>0}$ with the topology
whose base of open sets consists of the intervals $]a,b[ \,\subset
\RR_{>0}$ and the subsets of the form $]a, \infty] = \; ]a, \infty[ \,
\cup  \, \{\infty \}$. We denote the Euclidean norm on $\RR^n$ by $\| \cdot \|$
and the unit sphere in $\RR^n$ by $S_1$. 
\begin{proposition}  \label{l:star}
Endow $\RR^n$ with its standard smooth structure, let $\mu
\colon S_1 \ra \overline{\RR}_{>0}$ be a continuous function and let
\[
S \,=\, \left\{ x \in \RR^n \, \Big|\,  0 \le \| x \| < 
\mu \left( \frac{x}{\| x \|} \right) \right\}
\]
be the starlike domain associated with $\mu$.
Then $S$ is diffeomorphic to $\RR^n$.
\end{proposition}

\begin{remark} 
{\rm
The diffeomorphism guaranteed by Proposition \ref{l:star} may 
be chosen such that the rays emanating from the origin are preserved.
}
\end{remark}

\noindent
{\it Proof of Proposition \ref{l:star}.} 
If $\mu (S_1) = \{ \infty \}$, there is nothing to prove. 
In the case that $\mu$ is bounded, \text{Proposition \ref{l:star}} 
has been proved by Ozols \cite{O}. 
In the case that neither $\mu (S_1) = \{ \infty \}$ nor $\mu$ is bounded,
Ozols's proof readily extends to this situation. 
Using his notation, the
only modifications needed are: Require in addition that $r_0 <1$ and
that $\eps_1 <2$, and define continuous functions $\tilde{\mu}_i \colon
S_1 \ra \RR_{>0}$ by 
\[
\tilde{\mu}_i = \min \left\{ i ,\, \mu - \eps_i + \tfrac{\dd_i}{2} \right\}.
\]
With these minor adaptations the proof in \cite{O} applies word by word.
\proofend

In the following we shall use some basic Riemannian geometry.
We refer to \cite{KN} for basic notions and results in Riemannian geometry.
Consider an $n$-dimensional complete Riemannian manifold $(N, g)$.
We denote the cut locus of a point $p \in N$ by $C(p)$.
\begin{corollary}  \label{c:full}
The maximal normal neighbourhood $N \setminus C(p)$ of any point $p$ in
an $n$-dimensional complete Riemannian manifold $(N,g)$ is
diffeomorphic to $\RR^n$ endowed with its standard smooth structure. 
\end{corollary}

\proof
Fix $p \in N$. We identify the tangent space $(T_pN, g(p))$ with
Euclidean space $\RR^n$ by a (linear) isometry.
Let $\exp_p \colon \RR^n \ra N$ be the
exponential map at $p$ with respect to $g$, 
and let $S_1$ be the unit sphere in $\RR^n$.
We define the function $\mu \colon S_1 \ra \overline{\RR}_{>0}$ by
\begin{equation}  \label{d:mu}
\mu (x) \,=\, \inf \{ \:\! t>0 \mid \exp_p (tx) \in C(p) \:\! \} .
\end{equation}
Since the Riemannian metric $g$ is complete, the function $\mu$ is
continuous \cite[VIII, Theorem 7.3]{KN}.
Let $S \subset \RR^n$ be the starlike domain associated with $\mu$. 
In view of Proposition \ref{l:star} the set $S$ is diffeomorphic to $\RR^n$,
and in view of \cite[VIII, Theorem 7.4 \:\!(3)]{KN} we have
$\exp_p(S) = N \setminus C(p)$.
Therefore, $N \setminus C(p)$ is diffeomorphic to $\RR^n$.
\proofend

A main ingredient of our proof of Theorem \ref{ta:1}
are the following two special 
cases of a theorem of Greene and Shiohama \cite{GS}.
\begin{proposition} \label{l:GS}
(i) Assume that $\Omega_1$ is a volume form on the connected open subset
$U$ of $\RR^n$ such that 
$\Vol (U, \Omega_1) = \left| U \right| < \infty$.
Then there exists a diffeomorphism $\psi$ of $U$ such that $\psi^* \Omega_1 = \Omega_0$.

\s
(ii) Assume that $\Omega_1$ is a volume form on $\RR^n$ such that
$\Vol (\RR^n, \Omega_1) = \infty$.
Then there exists a diffeomorphism $\psi$ of $\RR^n$ such that 
$\psi^* \Omega_1 = \Omega_0$.
\end{proposition}

\m
\ni
{\bf End of the proof of Theorem \ref{ta:1}.}

\s
\ni
Let $U \subset \RR^n$ and $(M, \Omega)$ be as in Theorem \ref{ta:1}.
After enlarging $U$, if necessary, we can assume that 
$\left| U \right| = \Vol (M, \Omega)$.
We set $N = M \setminus \pp M$.
Then 
\begin{equation}  \label{e:vol:NM}
\left| U \right| \,=\, \Vol (M, \Omega) \,=\, \Vol (N, \Omega) .
\end{equation}
Since $N$ is a connected manifold without boundary, there exists 
a complete Riemannian metric $g$ on $N$.  
Indeed, according to a theorem of Whitney \cite{W}, 
$N$ can be embedded as a closed submanifold in some $\RR^m$. 
We can then take the induced Riemannian metric. A direct and elementary
proof of the existence of a complete Riemannian metric is given in
\cite{NO}.

Fix a point $p \in N$.
As in the proof of Corollary \ref{c:full} we identify $(T_pN, g(p))$
with $\RR^n$ and define the
function $\mu \colon S_1 \ra \overline{\RR}_{>0}$ as in \eqref{d:mu}.
Using polar coordinates on $\RR^n$ we see from Fubini's Theorem that the
set
\[
\widetilde{C}(p) \,=\, \{ \:\! \mu(x) x \mid x \in S_1 \:\! \}
\,\subset\, \RR^n
\]
has measure zero, and so 
$C(p) = \exp_p \big( \widetilde{C}(p) \big)$ also has measure zero 
(see \cite[VI, Corollary 1.14]{Bo}). 
It follows that
\begin{equation}  \label{e:vol:NC}
\Vol (N \setminus C(p), \Omega) \,=\, \Vol (N, \Omega).
\end{equation}
According to Corollary \ref{c:full} there exists a diffeomorphism
\[
\dd \colon \RR^n \,\ra\, N \setminus C(p) .
\]
After composing $\dd$ with a reflection of $\RR^n$, if necessary, we can 
assume that $\dd$ is orientation preserving.
In view of \eqref{e:vol:NM} and \eqref{e:vol:NC} we then have
\begin{equation}  \label{e:vol:Ud}
\left| U \right| \,=\, \Vol (\RR^n, \dd^* \Omega).
\end{equation}

\s
\ni
{\bf Case 1.}
$\left| U \right| < \infty$.

\s
\ni
Let $U_i$, $i = 1,2, \dots$, be the countably many components of $U$. 
Then $0 < \left| U_i \right| < \infty$ for each $i$. 
Given numbers $a$ and $b$ with $- \infty \le a < b \le \infty$ we
abbreviate the ``open strip''
\[
S_{a,b} = \{\:\! (x_1, \dots, x_n) \in \RR^n \mid a<x_1<b \:\!\} .
\]
In view of the identity \eqref{e:vol:Ud} we have
\[
\sum_{i \ge 1} \left| U_i \right| \,=\, \left| U \right| \,=\, \Vol
(\RR^n, \dd^* \Omega) .
\]
We can therefore inductively define $a_0 = - \infty$ and
$a_i \in \:]-\infty, \infty]$ by
\[
\Vol \left( S_{a_{i-1}, a_i}, \dd^* \Omega \right) \,=\, \left| U_i \right| .
\]
Abbreviating $S_i = S_{a_{i-1}, a_i}$ we then have
$\RR^n = \bigcup_{i \ge 1} \overline{S_i}$.

For each $i \ge 1$ we choose an orientation preserving diffeomorphism
$\tau_i \colon \RR^n \ra S_i$.
In view of Proposition \ref{l:smooth} we find a smooth embedding $\ss_i \colon U_i \ha \RR^n$ such that $\RR^n \setminus \ss_i (U_i)$ has measure zero.
After composing $\ss_i$ with a reflection of $\RR^n$, if necessary, we
can assume that $\ss_i$ is orientation preserving.
Using the definition of the volume, we can now conclude that
\[
\Vol (U_i, \ss_i^* \tau_i^* \dd^* \Omega) = 
\Vol ( \ss_i(U_i), \tau_i^* \dd^* \Omega) = 
\Vol ( \RR^n, \tau_i^* \dd^* \Omega) = 
\Vol (S_i, \dd^* \Omega) = \left| U_i \right| . 
\]
In view of Proposition \ref{l:GS}\:(i) we therefore find a diffeomorphism
$\psi_i$ of $U_i$ such that 
\begin{equation}  \label{e:vol:pst}
\psi_i^* \left( \ss_i^* \tau_i^* \dd^* \Omega \right) \,=\, \Omega_0 .
\end{equation}
We define $\ff_i \colon U_i \ha M$ to be the composition of
diffeomorphisms and smooth embeddings 
\[
U_i \,\xrightarrow{\psi_i}\, U_i 
    \,\xrightarrow{\ss_i}\, \RR^n 
    \,\xrightarrow{\tau_i}\, S_i \,\subset\, \RR^n 
    \,\xrightarrow{\dd} N \setminus C(p) \,\subset\, M .
\]
The identity \eqref{e:vol:pst} implies that $\ff_i^* \Omega = \Omega_0$.
The smooth embedding
\[
\ff = \coprod \ff_i \colon U = \coprod U_i \,\ha\, M
\]
therefore satisfies $\ff^* \Omega = \Omega_0$.

\m
\ni
{\bf Case 2.}
$\left| U \right| = \infty$.

\s
\ni
In view of \eqref{e:vol:Ud} we have $\Vol (\RR^n, \dd^* \Omega) = \infty$.
Proposition \ref{l:GS}\:(ii) shows that there exists a diffeomorphism $\psi$ of $\RR^n$ such that 
\begin{equation}  \label{e:psi2}
\psi^* \dd^* \Omega = \Omega_0 .
\end{equation}
We define $\ff \colon U \ha M$ to be the composition of inclusions 
and diffeomorphisms
\[
U \,\subset\, \RR^n \,\xrightarrow{\psi}\, \RR^n \,\xrightarrow{\dd}\, N 
\setminus C(p) \,\subset\, M .
\] 
The identity \eqref{e:psi2} implies that $\ff^* \Omega = \Omega_0$.
The proof of Theorem \ref{ta:1} is complete.
\proofend

\end{document}